\documentclass[11pt]{amsart}
\usepackage{epsfig}
\usepackage{amssymb}
\usepackage{amscd}
\usepackage{epsf}

\long\def\comment#1\endcomment{\relax}

\newcounter{subsubsubsection}
\newcounter{subsubsubsubsection}

\makeatletter
\@addtoreset{subsubsubsection}{subsubsection}
\@addtoreset{subsubsubsubsection}{subsubsubsection}

\newcommand{\subsubsubsection}[1]{\par\addtocounter{subsubsubsection}{1}\smallskip
{\thesubsubsection.\arabic{subsubsubsection}.\bfseries\ #1}}

\newcommand{\subsubsubsubsection}[1]{\par\addtocounter{subsubsubsubsection}{1}\smallskip
{\thesubsubsection.\arabic{subsubsubsection}.\arabic{subsubsubsubsection}\bfseries\ #1}}

\makeatother

\newcommand{\sevafig}[3]{\begin{figure}[h]\centerline{
 \epsfig{file=#1,width=#2,angle=#3}}
\bigskip\caption{}\end{figure}}

\newcommand{\one}{{\mathbf1}}

\renewcommand{\labelenumi}{(\roman{enumi})\,}

\newcommand{\U}{\mathcal U}

\renewcommand{\Im}{\matho{Im}\nolimits}

\newcommand{\PBW}{{PBW}}

\newcommand{\GL}{G\!L}
\newcommand{\simto}{\overset{\sim}{\to}}

\newcommand{\DR}{{\mathrm{DR}}}

\newcommand{\inv}{{\mathrm{inv}}}
\newcommand{\coinv}{{\mathrm{coinv}}}

\newcommand{\ndot}{\bullet}

\newcommand{\g}{{\mathfrak g}}
\newcommand{\m}{{\mathfrak m}}

\def\matho#1{\mathop{\mathrm{#1}}}

\newcommand*{\ato}[2]{{\genfrac{}{}{0pt}{}{#1}{#2}}}

\DeclareMathSizes{11.1}{10}{8}{6}

\newcommand{\Alt}{\matho{Alt}}

\newcommand{\Hom}{\matho{Hom}\nolimits}

\newcommand{\Star}{\mathrm{Star}}

\newcommand{\Conf}{\mathrm{Conf}}
\newcommand{\opp}{{\mathrm{opp}}}

\newcommand{\Tr}{{\matho{Tr}}}
\newcommand{\pr}{\matho{pr}}

\newcommand{\Id}{\matho{Id}\nolimits}
\newcommand{\ad}{\matho{ad}}
\newcommand{\Log}{\matho{Log}}

\newcommand{\strange}{{\mathrm{strange}}}
\newcommand{\poly}{{\mathrm{poly}}}
\newcommand{\Disk}{\mathrm{Disk}}

\newcommand{\HKR}{{\mathrm{HKR}}}

\newcommand{\Hoch}{\mathrm{Hoch}}

\newcommand{\C}{\mathbb C}
\newcommand{\R}{\mathbb R}
\newcommand{\Z}{\mathbb Z}
\newcommand{\D}{\mathcal D}
\renewcommand{\H}{\mathcal H}

\newtheorem*{theorem}{Theorem}

\newtheorem*{lemma}{Lemma}
\newtheorem*{corollary}{Corollary}
\newtheorem*{conjecture}{Conjecture}

\theoremstyle{remark}

\newtheorem{nexample}{Example}

\theoremstyle{definition}
\newtheorem*{defin}{Definition}

\author{Boris Shoikhet}
\title
{A proof of the Tsygan formality conjecture for~chains}
\date{2000}
\address{FIM, ETH-Zentrum, CH-8092 Z\"urich, SWITZERLAND}
\email{borya@mccme.ru}

\begin{document}

\sloppy

\maketitle

\begin{abstract}
We extend the Kontsevich formality $L_\infty$-morphism $\U\colon
T^\ndot_\poly(\R^d)\to\D^\ndot_\poly(\R^d)$ to an $L_\infty$-morphism
of $L_\infty$-modules over $T^\ndot_\poly(\R^d)$, $\hat \U\colon
C_\ndot(A,A)\to\Omega^\ndot(\R^d)$, $A=C^\infty(\R^d)$. The
construction of the map~$\hat \U$ is given in Kontsevich-type
integrals. The conjecture that such an $L_\infty$-morphism exists is
due to Boris Tsygan \cite{Ts}. As an application, we obtain an
explicit formula for isomorphism $A_*/[A_*,A_*]\simto A/\{A,A\}$
($A_*$~is the Kontsevich deformation quantization of the algebra~$A$
by a Poisson bivector field, and $\{{,}\}$ is the Poisson bracket). We
also formulate a conjecture extending the Kontsevich theorem on 
cup-products to this context. The conjecture implies a generalization
of the Duflo formula, and many other things.
\end{abstract}

\section{$L_\infty$-algebras and $L_\infty$-modules}

Here we recall basic definitions from~\cite{Ts} and construct an
$L_\infty$-module over $T^\ndot_\poly(\R^d)$ structure on the chain
Hochschild complex $C_\ndot(A,A)$, $A=C^\infty(\R^d)$.

\subsection{}
For a manifold~$M$, $T^\ndot_\poly(M)$ denotes the graded Lie algebra
of smooth polyvector fields on~$M$,  $\Omega^\ndot(M)$ denotes the
graded space of smooth differential forms on~$M$. The space
$T^\ndot_\poly(M)$ is graded as a Lie algebra, i.e.\
$T^i_\poly(M)=\{\text{$(i+1)$-polyvector fields}\}$, $i\ge-1$. The bracket is
the Schouten--Nijenhuis bracket, generalizing the usual Lie bracket of
vector fields. It is defined as follows:
\begin{multline}
[\xi_0\wedge\dots\wedge\xi_k,\eta_0\wedge\dots\wedge\eta_l]=\\
=\sum_{i=0}^k\sum_{j=0}^l(-1)^{i+j+k}[\xi_i,\eta_j]\wedge\xi_0\wedge\dots\wedge
\hat\xi_i\wedge\dots\wedge\xi_k\wedge\eta_0\wedge\dots\wedge\hat\eta_j\wedge\dots\wedge\eta_l
\end{multline}
for $k,l\ge0$ ($\{\xi_i\}$ and $\{\eta_j\}$ are vector fields, the
bracket~(1) does not depend on their choice and depends only on
polyvector fields $\xi=\xi_0\wedge\dots\wedge\xi_k$ and
$\eta=\eta_0\wedge\dots\wedge\eta_l$) and for $k\ge0$
\begin{equation}
[\xi_0\wedge\dots\wedge\xi_k,h]=\sum_{i=0}^k(-1)^i\xi_i(h)\xi_0\wedge\dots\wedge\hat\xi_i\wedge\dots\wedge\xi_k
\end{equation}
(here $h\in C^\infty(M)$ is a function).

It is clear that $[T^i_\poly,T^j_\poly]\subset T^{i+j}_\poly$. We will
consider $T^\ndot_\poly$ as a differential graded (dg) Lie algebra
equipped with zero differential.

The space $\Omega^\ndot(M)$ is $\Z_{\le0}$-graded:
$$
\Omega^i(M)=\{\text{$i$-differential forms}\},\quad \deg\Omega^i=-i.
$$

For a $k$-polyvector field $\gamma$ we denote by $i_\gamma$ the
natural contraction
$$
i_\gamma\colon\Omega^\ndot(M)\to\Omega^{\ndot-k}(M).
$$

The Lie derivative $L_\gamma$ is defined by the Cartan formula
\begin{equation}
L_\gamma=d_{\DR}\circ i_\gamma+i_\gamma\circ d_\DR=[d_\DR,i_\gamma]
\end{equation}
(we denote $[a,b]_+$ by $[a,b]$ for odd symbols $a,b$).

\begin{lemma}
$[L_{\gamma_1},L_{\gamma_2}]=L_{[\gamma_1,\gamma_2]}$ where
$[\gamma_1,\gamma_2]$ is the Schouten--Nijenhuis bracket \emph{(1),
(2)} of polyvector fields.
\end{lemma}

In this way we endow the graded space $\Omega^\ndot(M)$ with a structure
of a module over the graded Lie algebra $T^\ndot_\poly(M)$. Note that
for $\gamma\in T^k_\poly(M)$ one has
$L_\gamma\Omega^\ndot(M)\subset\Omega^{\ndot-k}(M)$ and
$\deg\Omega^{\ndot-k}=\deg\Omega^\ndot+k$ because we consider
$\Omega^\ndot$ to be $\Z_{\le0}$ graded.

\subsection{}
We denote by $C^\ndot(A,A)$ the cohomological Hochschild complex of an
associative algebra~$A$. By definition, $C^k(A,A)=\Hom_\C(A^{\otimes
k},A)$, $k\ge0$. The differential $d_\Hoch\colon C^k(A,A)\to
C^{k+1}(A,A)$ is defined as follows:
\begin{multline}
(d\Psi)(a_1\otimes\dots\otimes
a_{k+1})=a_1\cdot\Psi(a_2\otimes\dots\otimes a_{k+1})-\\
-\Psi(a_1\cdot a_2\otimes a_3\otimes\dots\otimes
a_{k+1})+\dots\pm\Psi(a_1\otimes\dots\otimes a_k)\cdot a_{k+1}.
\end{multline}
By definition,
$\D^\ndot_\poly(M)=C^\ndot_\poly(C^\infty(M),C^\infty(M)){[1]}$ is
shifted by~$1$ subcomplex of $C^\ndot(A,A)$, $A=C^\infty(M)$,
consisting of polydifferential operators. Thus,
$\D^k_\poly(M)\subset\Hom_\C(A^{\otimes k+1},A)$, $k\ge-1$.

The Gerstenhaber bracket on $C^\ndot(A,A)[1]$ is defined as follows:
for $\phi_1\in C^{k_1+1}(A,A)$, $\phi_2\in C^{k_2+1}(A,A)$
\begin{equation}
[\phi_1,\phi_2]=\phi_1\circ\phi_2-(-1)^{k_1k_2}\phi_2\circ\phi_1
\end{equation}
where
\begin{multline}
(\phi_1\circ\phi_2)(a_0\otimes\dots\otimes a_{k_1+k_2})=\\
=\sum_{i=0}^{k_1}(-1)^{ik_2}\phi_1(a_0\otimes\dots\otimes
a_{i-1}\otimes\phi_2(a_i\otimes\dots\otimes a_{i+k_2})\otimes
a_{i+k_2+1}\otimes\dots\otimes a_{k_1+k_2}).
\end{multline}

\begin{lemma}
Formulas \emph{(4), (5), (6)} define a dg Lie algebra structure on
$C^\ndot(A,A)[1]$.
\end{lemma}

Note that
\begin{equation}
d_\Hoch(\Psi)=[m,\Psi]
\end{equation}
where $m\colon A^{\otimes2}\to A$ is the product, the associativity is
equivalent to $[m,m]=0$, and $\ad m$ defines a differential.

\subsubsection{}

\begin{lemma}[Hochschild--Kostant--Rosenberg]
\leavevmode

\begin{enumerate}
\item $H^i(\D^\ndot_\poly(M))=T^i_\poly(M)$, and the map
$\varphi_{\HKR}\colon T^\ndot_\poly\to\D^\ndot_\poly$,
\begin{equation}
\varphi_{\HKR}(\xi_1\wedge\dots\wedge\xi_k)(f_1\otimes\dots\otimes
f_k)=\frac1{k!}\Alt_{\xi_1,\dots,\xi_k}\prod_{i=1}^k\xi_i(f_i)
\end{equation}
is a quasi-isomorphism of complexes \emph(in particular,
$d_\Hoch|_{\Im\varphi_\HKR}=0$\emph{);}
\item the induced map on cohomology $[\varphi_\HKR]\colon
T^\ndot_\poly(M)\to H^\ndot(\D^\ndot_\poly(M))$ 
is a Lie algebras isomorphism.
\end{enumerate}
\end{lemma}

Formality theorem of M.~Kontsevich~\cite{K} states that
$T^\ndot_\poly(M)$ and $\D^\ndot_\poly(M)$ are 
 quasi-isomorphic as dg
Lie algebras, in the sense of the derived categories, or
$L_\infty$-quasi-isomorphic.

\subsection{}
We denote by $C_\ndot(A,A)$ the homological (chain) Hochschild complex
of an associative algebra~$A$. By definition, $C_k(A,A)=A\otimes
A^{\otimes k}$, and the differential $b\colon C_k(A,A)\to
C_{k-1}(A,A)$ is defined as follows:
\begin{multline}
b(a_0\otimes\dots\otimes a_k)=\\
=a_0a_1\otimes a_2\otimes\dots\otimes
a_k-a_0\otimes a_1a_2\otimes\dots\otimes a_k+\dots\pm a_ka_0\otimes
a_1\otimes\dots\otimes a_{k-1}.
\end{multline}
We consider $C_\ndot(A,A)$ to be $\Z_{\le0}$-graded, $\deg(A\otimes
A^{\otimes k})=-k$, $k\ge0$.

\subsubsection{}
\begin{theorem}
\leavevmode\par
\begin{enumerate}
\item $H_i(C_\ndot(A,A))=\Omega^i(M)$, $A=C^\infty(M)$\emph;
\item the map
\begin{equation}
\begin{gathered}
\mu\colon C_\ndot(A,A)\to\Omega^\ndot(M),\\
\mu(a_0\otimes a_1\otimes\dots\otimes
a_k)=\frac1{k!}a_0\,da_1\wedge\dots\wedge da_k
\end{gathered}
\end{equation}
is a quasi-isomorphism of the complexes \emph($\Omega^\ndot(M)$ is
equipped with zero differential\emph).
\end{enumerate}
\end{theorem}

\comment

\subsubsection{}
We will need the following variation of the Theorem 1.3.1(ii) which
seems to be new.

\begin{theorem}
The map $\mu_!\colon C_\ndot(A,A)\to\Omega^\ndot(M)$ defined as
\begin{multline*}
\mu_!(a_0\otimes a_1\otimes\dots\otimes
a_k)=\frac1{k!}(a_0\,da_1\wedge\dots\wedge da_k+(-1)^ka_1\,da_2\wedge
da_3\wedge\dots\wedge da_k\wedge da_0+\\
\begin{gathered}+(-1)^k\dots\pm a_k\,da_0\wedge
da_1\wedge\dots\wedge da_{k-1})\\
\text{\emph(all cyclic permutations\emph)}
\end{gathered}
\end{multline*}
is quasi-isomorphism of the complexes \emph($\Omega^\ndot(M)$ is
equipped with zero differential\emph).
\end{theorem}

We will use this result in the sequel, but we omit the proof.

\endcomment

\subsection{}

In Section 1.1 we have defined operators~$L_\gamma$, $\gamma\in
T^\ndot_\poly(M)$, acting on $\Omega^\ndot(M)$. We know from Theorems~1.2.1(i), 1.3.1(i) that
$H^\ndot(\D^\ndot_\poly(M))=T^\ndot_\poly(M)$,
$H_\ndot(C_\ndot(A,A))=\nobreak\Omega^\ndot(M)$, $A=C^\infty(M)$. In this
section we define the operators $L_\gamma$ ``on the level of
complexes'', i.e., we define operators~$L_\Psi$, $\Psi\in
C^\ndot(A,A)[1]$, acting on $C_\ndot(A,A)$. The operator~$L_\Psi$,
$\Psi\in \Hom(A^{\otimes k},A)$, is defined as follows:
\begin{multline}
L_\Psi(a_0\otimes\dots\otimes a_n)=
\sum_{i=0}^{n-k}(-1)^{(k-1)(i+1)}a_0\otimes\dots\otimes
a_i\otimes\Psi(a_{i+1}\otimes\dots\otimes a_{i+k})\otimes\dots\otimes
a_n+\\
+\sum_{j=n-k}^n(-1)^{n(j+1)}\Psi(a_{j+1}\otimes\dots\otimes
a_0\otimes\dots)\otimes a_{k+j-n}\otimes\dots\otimes a_j.
\end{multline}

Note that
\begin{equation}
L_m(a_0\otimes\dots\otimes a_n)=b(a_0\otimes\dots\otimes a_n)
\end{equation}
(here $m\colon A^{\otimes2}\to A$ is the multiplication, and $b$ is
the chain Hochschild differential, see~(9)).

\subsubsection{}
\begin{lemma}
$$
[L_{\Psi_1},L_{\Psi_2}]=L_{[\Psi_1,\Psi_2]}
$$
where $[\Psi_1,\Psi_2]$ is the Gerstenhaber bracket of the cochains
$\Psi_1$ and~$\Psi_2$.
\end{lemma}

In particular, let $\Psi_1=m\colon A^{\otimes2}\to A$ be the
multiplication; we know that $L_m=b$, and we have:
\begin{equation}
[b, L_\Psi]=L_{d_\Hoch\Psi}
\end{equation}
(see~(7)).

Lemma 1.4.1 allows us to equip the space $C_\ndot(A,A)$ with a
structure of the module over the dg Lie algebra $C^\ndot(A,A)[1]$. We
explain below how, using an $L_\infty$-map $\U\colon
T^\ndot_\poly\to\D^\ndot_\poly$, to equip the space $C_\ndot(A,A)$
with an $L_\infty$-module structure over~$T^\ndot_\poly$.

\subsection{}

In this Subsection we recall the basic definitions of homotopical
algebra. We will do it ``on the level of formulas'', leaving the
language of formal $Q$-manifolds~\cite{K}, because it is more
convenient for our needs.

\subsubsection{}
An $L_\infty$-algebra is a $\Z$-graded vector space~$\g$ with a
collection of maps:
\begin{align*}
&Q_1\colon\g\to\g[1]\\
&Q_2\colon\Lambda^2\g\to\g\\
&Q_3\colon\Lambda^3\g\to\g[-1]\\
&\hbox to 3cm{\dotfill}
\end{align*}
satisfying the relations
\begin{equation}
\sum_{\ato{i_1<\dots<i_p,\ j_1<\dots<j_q}{p+q=k}}\pm
Q_{q+1}(Q_p(x_{i_1}\wedge\dots\wedge x_{i_p})\wedge
x_{j_1}\wedge\dots\wedge x_{j_p})=0
\end{equation}
for each $k\ge2$ and homogenous $\{x_s\}$.

The first relation, for $k=2$, is
$$
Q^2_1=0.
$$

The second, $k=3$, is that the product $Q_2$ is compatible with the
differential~$Q_1$, i.e.
$$
Q_1(Q_2(x\wedge y))=Q_2(Q_1x\wedge y)\pm Q_2(x\wedge Q_1y).
$$

The third is that the skew-symmetric product $Q_2$ obeys the Jacobi
identity modulo~$Q_3$. The case when $Q_3=Q_4=\dots=0$ is the case of
dg Lie algebras.

\subsubsection{}
An $L_\infty$-morphism $\U\colon\g_1\to\g_2$ between two
$L_\infty$-algebras is a collection of maps:
\begin{align*}
&\U_1\colon\g_1\to\g_2\\
&\U_2\colon\Lambda^2\g_1\to\g_2[-1]\\
&\U_3\colon\Lambda^3\g_1\to\g_2[-2]\\
&\hbox to 3cm{\dotfill}
\end{align*}
obeying the identities
\begin{multline}
\sum\pm \U_{q+1}(Q_p(x_{i_1}\wedge\dots\wedge x_{i_p})\wedge
x_{j_1}\wedge\dots\wedge x_{j_q})=\\
\begin{gathered}
=\sum\pm\frac1{k!}Q_k(\U_{n_1}(x_{i_{11}}\wedge\dots\wedge
x_{i_{1n_1}})\wedge\dots\wedge \U_{n_k}(x_{i_{k1}}\wedge\dots\wedge
x_{i_{kn_k}}))\\
(n_1+\dots+n_k=p+q).
\end{gathered}
\end{multline}

In the simplest case, when $\g_1$ and $\g_2$ are dg Lie algebras,
$Q_1=d$, $Q_2=[{,}]$, (15)~is:
\begin{multline}
d\U_n(x_1\wedge\dots\wedge x_n)+\sum_{i=1}^n\pm
\U_n(x_1\wedge\dots\wedge dx_i\wedge\dots\wedge x_n)=\\
=\frac12\sum_{\ato{k,l\ge1}
{k+l=n}}\frac1{k!l!}\sum_{\sigma\in\Sigma_n}\pm[\U_k(x_{\sigma_1}\wedge\dots\wedge
x_{\sigma_k}),\U_l(x_{\sigma_{k+1}}\wedge\dots\wedge
x_{\sigma_n})]+\\
+\sum_{i<j}\pm \U_{n-1}([x_i,x_j]\wedge
x_1\wedge\dots\wedge\hat x_i\wedge\dots\wedge\hat x_j\wedge\dots\wedge
x_n).
\end{multline}

The simplest cases are:

\begin{description}
\item[$n=1$] $d\U_1(x)=\U_1(dx)$, i.e.\ $\U_1$ is a map of 
complexes; an $L_\infty$-morphism is called an
$L_\infty$-quasi-isomorphism if $\U_1$ is a quasi-isomorphism of 
complexes;
\item[$n=2$] $\U_1$ is a map of dg Lie algebras modulo~$\U_2$; $\U_1$ is
a map of graded Lie algebras on the level of cohomology.
\end{description}

\noindent
The connection of the notion of an $L_\infty$-morphism with the
classical homological algebra is that if two dg Lie algebras $\g_1$
and $\g_2$ are quasi-isomorphic in the sense of derived cathegeries,
that is there exists a hat
$\g_1\overset{\varphi_1}\to\g_3\overset{\varphi_2}\leftarrow\g_2$ for some
dg Lie algebra~$\g_3$, and both maps $\varphi_1$ and $\varphi_2$ are maps of
dg Lie algebras and quasi-isomorphisms of the complexes, they are $L_\infty$-quasi-isomorphic. In the $L_\infty$-side
we do not construct an extra dg Lie algebra $\g_3$ but we construct an
infinitely many higher ``Taylor components'' of the map~$\U_1$.

\subsubsection{}
An $L_\infty$-module $M$ over an $L_\infty$-algebra $\g$ is a graded
vector space $\g$ with a collection of maps
$$
\varphi_k\colon\Lambda^k\g\otimes M\to M[1-k],\quad k\ge0
$$
satisfying for each $k\ge0$ the equation
\begin{multline}
\sum_{p+q=k}\sum_{\ato{i_1<\dots<i_p}{
j_1<\dots<i<j_q}}\pm\phi_p(x_{i_1}\wedge\dots\wedge
x_{i_p}\otimes\phi_q(x_{j_1}\wedge\dots\wedge x_{j_q}\otimes
m))+\\
+\sum_{p+q=k}\sum_{\ato{i_1<\dots<i_p}{
j_1<\dots<j_q}}\pm\phi_{q+1}(Q_p(x_{i_1}\wedge\dots\wedge
x_{i_p})\wedge x_{j_1}\wedge\dots\wedge x_{j_q}\otimes m)=0.
\end{multline}

For $k=0$ (17) gives $\phi^2_0=0$, i.e.\ $\phi_0$ is a differential on
th graded space~$M$.

For $k=1$ we obtain that the map $\phi_1\colon\g\otimes M\to M$ is a
map of the complexes:
$$
\phi_0(\phi_1(x\otimes m))=\phi_1(dx\otimes
m)\pm\phi_1(x\otimes\phi_0m).
$$

In the case when $\g$ is a dg Lie algebra we obtain that the map
$\phi_1$ defines a $\g$-module structure on $M$ modulo~$\phi_2$, etc.

\subsubsection{}

A morphism $\varphi$ of two $L_\infty$-modules $M,N$ over an
$L_\infty$-algebra~$\g$ is a collection of maps
$$
\varphi_k\colon\Lambda^k\g\otimes M\to N[-k],\quad k\ge0
$$
satisfying
\begin{multline}
\sum\pm\varphi_{q+1}(Q_p(x_{i_1}\wedge\dots\wedge x_{i_p})\wedge
x_{j_1}\wedge\dots\wedge x_{j_q}\otimes
m)+\\
+\sum\pm\varphi_{p+1}(x_{i_1}\wedge\dots\wedge
x_{i_p}\wedge\phi_q(x_{j_1}\wedge\dots\wedge x_{j_q}\otimes m))=\\
=\sum\pm\phi_{p+1}(x_{i_1}\wedge\dots\wedge
x_{i_p}\wedge\varphi_q(x_{j_1}\wedge\dots\wedge x_{j_q}\otimes m))
\end{multline}
for $p+q=k\ge0$.

For $k=p+q=0$ we obtain that the map $\varphi_0$ is a map of the
complexes, etc.

\subsection{An $L_\infty$-module over $T^\ndot_\poly(\R^d)$ structure
on $C_\ndot(A,A)$}

M.~Kontsevich constructed in~\cite{K} an
$L_\infty$-(quasi-iso)morphism $\U\colon
T^\ndot_\poly(\R^d)\to\D^\ndot_\poly(\R^d)$. Let $\U_1,\U_2,\U_3,\dots$
be its Taylor components. Define a set of maps
$$
\phi_k\colon\Lambda^kT^\ndot_\poly(\R^d)\otimes C_\ndot(A,A)\to
C_\ndot(A,A)[1-k]
$$
by the formulas:
\begin{equation}
\left.
\begin{aligned}\relax
&\phi_0=b\ \text{(see (9))},\\
&\phi_k(\gamma_1\wedge\dots\wedge\gamma_k\otimes\omega)=L_{\Psi_k}\omega\\
&\ \text{with $\Psi_k=\U_k(\gamma_1\wedge\dots\wedge\gamma_k)$ for
$k\ge1$}
\end{aligned}
\right\}
\end{equation}

\begin{lemma}
In this way we have defined an $L_\infty$-module structure on
$C_\ndot(A,A)$ over dg Lie algebra $T^\ndot_\poly(\R^d)$.
\end{lemma}

\begin{proof}
We have to prove (17) for $\{\phi_k\}$.

For $k=0$ it is just $b^2=0$,

for $k=1$ it is
\begin{equation}
b\circ L_{\U_1(\gamma)}\omega=\pm L_{\U_1(\gamma)}\circ b\omega
\end{equation}
But we have from (13):
$$
[b,L_{\U_1(\gamma)}]\omega=L_{d_\Hoch \U_1(\gamma)}\omega
$$
and $d_\Hoch \U_1(\gamma)=\U_1(d_{T_\poly}\gamma)\equiv0$ for any
$L_\infty$-morphism $\U\colon T^\ndot_\poly\to\D^\ndot_\poly$. 

The general case is analogous. The statement of Lemma is true for any
$L_\infty$-morphism $\U\colon T^\ndot_\poly\to\D^\ndot_\poly$, not only for
the Kontsevich's one.
\end{proof}

\subsubsection{}
Now we have two modules over $T^\ndot_\poly(\R^d)$: the usual module
$\Omega^\ndot(\R^d)$, and the $L_\infty$-module $C_\ndot(A,A)$. We
want to construct an $L_\infty$-morphism
$$
\hat \U\colon C_\ndot(A,A)\to\Omega^\ndot(\R^d)
$$
between them.

It means that we search for maps
$$
\hat \U_k\colon\Lambda^kT^\ndot_\poly(\R^d)\otimes
C_\ndot(A,A)\to\Omega^\ndot(\R^d)[-k],\quad k\ge0
$$
such that for $k\ge-2$
\begin{multline}
\sum\pm\hat
\U_{k+1}([\gamma_{i_1},\gamma_{i_2}]\wedge\gamma_{j_1}\wedge\dots\wedge\gamma_{j_k}\otimes\omega)+
\sum_{\ato{p+q=k+2,}{
q\ge0}}\pm\hat
\U_{p+1}(\gamma_{i_1}\wedge\dots\wedge\gamma_{i_p}\otimes
L_{\U_q(\gamma_{j_1}\wedge\dots\wedge\gamma_{j_q})}\omega)+\\
+\sum\pm L_{\gamma_i}(\hat
\U_{k+1}(\gamma_{j_1}\wedge\dots\wedge\gamma_{j_{k+1}}\otimes\omega))=0
\quad(L_{\U_0}=b).
\end{multline}

We set:
$$
\hat \U_0(\omega)=\mu\omega\quad\text{(see Theorem 1.3.1)}.
$$

For $k=-2$ (21) is just the statement that $\mu$ is a map of the
complexes.

In the rest of this paper we construct such an
$L_\infty$-morphism~$\hat \U$, where $\U$ is the Kontsevich formality
morphism.

\section{Construction of the morphism $\hat \U$}

\subsection{The Kontsevich formality morphism $\U$}

Here we recall, very briefly, the construction~\cite{K} of the
Kontsevich formality morphism~$\U$.

We need to construct maps
$$
\U_k\colon\Lambda^kT^\ndot_\poly(\R^d)\to\D^\ndot_\poly(\R^d)[1-k].
$$

The formula for $\U_k$ is organized as a sum over the admissible
graphs~$\Gamma$. We cite from~\cite{K} the definition of an admissible
graph~$\Gamma$.

\subsubsection{}
\begin{defin}
Admissible graph $\Gamma$ is an oriented graph with labels such that
{\renewcommand{\labelenumi}{\arabic{enumi})}\begin{enumerate}
\item the set of vertices $V_\Gamma$ is
$\{1,\dots,n\}\sqcup\{\bar1,\dots,\bar m\}$, $n,m\in\Z_{\ge0}$,
$2n+2+m\ge0$; vertices from the set $\{1,\dots,n\}$ are called
vertices of the first type, vertices from $\{\bar1,\dots,\bar m\}$ are
called vertices of the second type,
\item every edge $(v_1,v_2)\in E_\Gamma$ starts at a vertex of the
first type,
\item there are no loops, i.e.\ no edges of the type $(v,v)$,
\item for every vertex $k\in\{1,\dots,n\}$ of the first type, the set
of edges
$$
\Star(k)=\bigl\{(v_1,v_2)\in E_\Gamma\mid v_1=k\bigr\}
$$
starting from~$k$, is labeled by symbols
$(e^1_k,\dots,e^{\#\Star(k)}_k)$.
\end{enumerate}}%
\end{defin}

\subsubsection{}
For any admissible graph $\Gamma$ with $n$ vertices of the first type
and $m$ vertices of the second type we define a map
$\U_\Gamma\colon\gamma_1\otimes\gamma_2\otimes\dots\otimes\gamma_n\mapsto\{A^{\otimes
m}\to A\}$, $A=C^\infty(\R^d)$, and $\gamma_i$ is
$(\#\Star(i))$-polyvector field. The function
$\Phi=\U_\Gamma(\gamma_1\otimes\dots\otimes\gamma_n)(f_1\otimes\dots\otimes
f_m)$ is defined as follows. It is the sum over all configurations of
indices running from $1$ to~$d$, labeled by~$E_\Gamma$:
$$
\Phi=\sum_{I\colon E_\Gamma\to\{1,\dots,d\}}\Phi_I
$$
where $\Phi_I$ is the product over all $n+m$ vertices of $\Gamma$ of
certain partial derivatives of functions $f_i$ and of coefficients
of~$\gamma_i$.

Namely, with each vertex $i$, $1\le i\le n$ of the first type we
associate function $\Psi_i$ an $\R^d$ which is
$$
\Psi_i=\langle\gamma_i,dx^{I(e_i^1)}\otimes\dots\otimes
dx^{I(e^{\#\Star(i)}_i)}\rangle.
$$
For each vertex $\bar j$ of the second type the associated function
$\Psi_{\bar j}$ is defined as~$f_j$.

In the next step we put into each vertex $v$ instead of
function~$\Psi_v$, its partial derivative
$$
\tilde\Psi_v=\left(\prod_{e\in E_\Gamma,\
e=(*,v)}\partial_{I(e)}\right)\Psi_v.
$$
Then, $\Phi_I=\prod_{v\in V_\Gamma}\tilde\Psi_v$, and
$\Phi=\sum_I\Phi_I$.

This is just a usual construction of $\GL_d$-invariants from graphs.

\subsubsection{}
Till now, the number of  edges of $\Gamma$ with $n$ vertices of the
first type and $m$ vertices of the second type was not  fixed. We
claim, that it is uniquely defined by the request that $\U_\Gamma$ is a
map
$$
\U_\Gamma\colon\Lambda^nT^\ndot_\poly(\R^d)\to\D^\ndot_\poly(\R^d)[1-n],
$$
i.e.\ by the grading.

As an element of $\D^\ndot_\poly(\R^d)$, $\U_\Gamma$ has grading
$m-1$. This grading should be equal to
$\deg\gamma_1+\dots+\deg\gamma_n+1-n$. We have:
$$
\deg\gamma_i=\#\Star(i)-1,
$$
and
$$
\#E_\Gamma=\sum_{i\in\{1,\dots,n\}}\#\Star(i),
$$
because any edge starts at a vertex of the first type by Definition
2.1.1.

Therefore,
$$
\#E_\Gamma=\deg\gamma_1+\dots+\deg\gamma_n+n.
$$

But
$$
m-1=\deg\gamma_1+\dots+\deg\gamma_n+1-n,
$$
therefore
$$
m-1=\#E_\Gamma-n+1-n,
$$
and
$$
\#E_\Gamma=2n+m-2.
$$

\subsubsection{}

Now we search a formula for $\U_n$ in the form
\begin{equation}
\U_n=\sum_\Gamma W_\Gamma\times \U_\Gamma,
\end{equation}
where $\Gamma$ has $n$ vertices of the first type, and $W_\Gamma\in\C$
is a number. We want to define $W_\Gamma$ as an integral of a form
$\Omega_\Gamma$ of the top degree over a configuration space. Any
edge of~$\Gamma$ will define a $1$-form on this configuration space,
and the form $\Omega_\Gamma$ is a wedge product of these $1$-forms (in
order corresponded to the labeling of the graph). Therefore, the
number $\#E_\Gamma$ of the edges should be equal to dimension of the
configuration space. Therefore, the configuration space should be a quotient by a
$2$-dimensional group of $n$ ``$2$-dimensional points'' and $m$
``$1$-dimensional'' points.

\subsubsection{}
Here we define these configuration spaces. Let $\H$ be the open
complex upper half-plane.

First, denote by $\Conf_{n,m}$ the space
\begin{multline*}
\Conf_{n,m}=\bigl\{p_1,\dots,p_n;q_1,\dots,q_m\mid p_i\in\H\
\text{for all}\ i=1,\dots, n,\\ 
q_j\in\R=\partial\bar\H\ \text{for
all}\ j=1,\dots,m;\
p_{i_1}\ne p_{i_2}\ \text{for}\ i_1\ne i_2,\ q_{j_1}\ne q_{j_2}\
\text{for}\ j_1\ne j_2\bigr\},
\end{multline*}
and
$$
\Conf_n=\bigl\{p_1,\dots,p_n\mid p_i\in\C, p_{i_1}\ne p_{i_2}\ \text{for}\ i_1\ne i_2\bigr\}.
$$
For $2n+m\ge2$ the group
$$
G^{(1)}=\bigl\{z\mapsto az+b\mid a\in\R_{>0},\ b\in\R\bigr\}
$$
acts freely on the space $\Conf_{n,m}$; therefore, the quotient-space
$$
C_{n,m}=\Conf_{n,m}/G^{(1)}
$$
is a smooth manifold of dimension $2n+m-2$.
Analogously,
$$
C_n=\Conf_n/G^{(2)}
$$
where
$$
G^{(2)}=\bigl\{z\mapsto az+b\mid a\in\R_{>0},\ b\in\C\bigr\}.
$$
It is a smooth manifold for $n\ge2$.

M.~Kontsevich constructed in~\cite{K}, Section~6, a compactification
of these spaces. The differential forms, defined below, can be
prolonged to regular differential forms on these
compactifications. Therefore, the integrals (defined below) converge.

We will not discuss these compactifications here. We will consider
this question for other spaces, suitable for the construction of
map~$\hat \U$, in Section~2.2.

\subsubsection{}
Now, for a graph $\Gamma$ with $n$ vertices of the first type and $m$
vertices of the second type, we are going to define a form
$\Omega_\Gamma$ of the top degree on~$C_{n,m}$, and then set
$W_\Gamma={\sim}\int_{C_{n,m}}\Omega_\Gamma$ (up to some coefficient
which depends on $n,m$ and does not depend on~$\Gamma$).

For any two points $p,q\in\bar\H$, denote by $l(p,q)$ the geodesic in
the Poincar\'e metric on~$\H$, passing through the points $p$
and~$q$. It is just a (part of) half-circle, passing through $p$
and~$q$, and orthogonal to $\R=\partial\bar\H$. Denote by
$l(p,\infty)$ the vertical line passing through the point~$p$. Denote
by $\phi^h(p,q)\in\R/2\pi\Z$ the angle between $l(p,q)$ and
$l(p,\infty)$, from $l(p,\infty)$ to $l(p,q)$ counterclockwise. The
formula for $\phi^h(p,q)$ is
\begin{equation}
\phi^h(p,q)=\frac1{2i}\Log\left(\frac{(q-p)(\bar q-p)}{(q-\bar p)(\bar
q-\bar p)}\right).
\end{equation}
The $1$-form $d\phi^h(p,q)$ is well-defined.

Now we set that the vertices $\{1,\dots,n\}$ of the first type of the
graph $\Gamma$ are the points $p_1,\dots,p_n$ on $\Conf_{n,m}$, and
vertices $\{\bar1,\dots,\bar m\}$ are the points $q_1,\dots,q_m$ on
$\Conf_{n,m}$. In other words, the graph ``is placed'' on $\H\sqcup\R$
in all possible ways. Then each edge $e=(p,q)$ of~$\Gamma$ defines a
map from $\Conf_{n,m}$ to $\Conf_{2,0}$ or to $\Conf_{1,1}$. On both spaces
$\Conf_{2,0}$ and $\Conf_{1,1}$ we have
constructed the $1$-form~$\phi^h$. Denote by $\varphi_e$ its pull-back to
$\Conf_{n,m}$. It is clear that $\varphi_e$ is $G^{(1)}$-invariant
because $G^{(1)}$ is the group of transformations in the Poincar\'e
metric on $\H$ preserving~$\{\infty\}$. Thus, $\varphi_e$ is a
well-defined $1$-form on $C_{n,m}=\Conf_{n,m}/G^{(1)}$.

We set:
\begin{equation}
W_\Gamma=\prod_{k=1}^n\frac1{(\#\Star(k))!}\cdot\frac1{(2\pi)^{2n+m-2}}\cdot\int_{C_{n,m}^+}\bigwedge\limits_{e\in
E_\Gamma}\varphi_e.
\end{equation}
Here $C^+_{n,m}$ is the connected component of $C_{n,m}$ where
$q_1<\dots<q_m$.

\begin{theorem}[M.~Kontsevich,\cite{K}]
The maps $\U_k$, $k\ge1$, defined by formulas \emph{(22), (24)}, are
the Taylor components of an $L_\infty$-morphism $\U\colon
T^\ndot_\poly(\R^d)\to\D^\ndot_\poly(\R^d)$. That is they satisfy the
identity~\emph{(16)}.
\end{theorem}

\subsubsection{}
One of the crucial points is that $G^{(1)}$ is the full group of
symmetries for the $1$-form $d\phi^h$. This is why the integrals (24)
a priori do not vanish (and they actually do not).

\subsection{Construction of the map $\hat \U$}

\subsubsection{Configuration spaces $D_{\one,n,m}$}

Let $D^2$ be the open disk $D^2=\bigl\{z\in\C\mid|z|<1\bigr\}$. We
denote by $\one$ the center of this disk $\one=\{z=0\}$. (It's a
notation bad from any point of view). The space $\Disk_{\one,n,m}$ is 
\begin{multline*}
\Disk_{\one,n,m}=\bigl\{p_1,\dots,p_n\in D^2,\ q_1,\dots,q_m\in
S^1=\partial\bar{D^2},\\
p_i\ne p_j\ \text{for}\ i\ne j,\ p_i\ne\one\ \text{for any}\ i,\ q_i\ne
q_j\ \text{for}\ i\ne j\bigr\}.
\end{multline*}
For $2n+m\ge1$ the group of rotations
$$
G=\bigl\{z\mapsto e^{i\theta}z,\ \theta\in\R/2\pi\Z\bigr\}
$$
acts freely on $\Disk_{\one,n,m}$, and we define the manifold
$D_{\one,n,m}$ as
\begin{equation}
D_{\one,n,m}=\Disk_{\one,n,m}/G.
\end{equation}
It has dimension $2n+m-1$. Note that $G\one=\one$.

Now we are going to describe a compactification $\bar\D_{\one,n,m}$
which will be used in the sequel.

The idea arise from~\cite{K}, Section~5; we just describe all strata
of codimension~$1$.

\begin{list}{}{}
\item[S1)] some points $p_{i_1},\dots,p_{i_k}$ of the first type move
close to each other and far from~$\one$, $k\ge2$. The corresponding
boundary stratum is $C_k\times D_{\one,n-k+1,m}$;
\item[S2)] some points $p_{i_1},\dots,p_{i_k}$ of the first type,
$k\ge1$, move close to each other and close to~$\one$. To describe
this stratum denote by $D_k$ the following space:
$$
D_k=\bigl\{p_1,\dots,p_k\in\C,\ p_i\ne p_j\ \text{for}\ i\ne
j\bigr\}\big/\bigl\{z\mapsto az,\ a\in\R_{>0}\bigr\}.
$$
For $k\ge1$ it is a manifold of dimension $k-1$. The boundary stratum
in the case~S2)  is $D_k\times D_{\one, n-k,m}$;
\item[S3)] some points $p_{i_1},\dots,p_{i_k}$ of the first type and
some points $q_{j_1},\dots,q_{j_l}$ of the second type move close to
each other (and close to $S^1=\partial\bar{D^2}$), $2k+l\ge2$. In this
case the boundary stratum is
$$
C_{k,l}\times D_{\one,n-k,m-l+1}.
$$
\end{list}

One can check that in all three cases the described strata have
codimension~$1$.

In the general case of arbitrary codimension, we have several groups of
points moving close to each other; this gives strata of
codimension~$>1$. Then ``we are looking through a magnifying glass''
to each group and find that inside it there are several groups of
points, moving close to each other, an so on. All the boundary strata
can be described by labeled and colored trees.

For us it is principal that the differential $1$-forms on
$D_{\one,n,m}$ (analogs of $d\phi^h$) constructed below can be
prolonged to regular differential $1$-forms on the
compactification $\bar D_{\one,n,m}$. The proof will be
clear after the definition of these $1$-forms.

\subsubsection{Admissible graphs}

An admissible graph is the same that in \cite{K} (see Section~2.1.1) but
here we have a marked vertex~$\one$, and \ $\one$ is not the end-point for any edge. There are
vertex~$\one$, vertices of the first type $\{1,\dots,n\}$, and
vertices of the second type $\{\bar1,\dots,\bar m\}$.

\subsubsection{Differential forms}

There are two types of edges: connecting $\one\in D^2$ with an other
point, and all others edges. We will define a differential $1$-form
$\varphi_e$ separately in these two cases.

\begin{list}{}{}
\item[C1)] $e$ is an usual edge, i.e.
$$
e=(p,q),\ p,q\in \bar{D^2},\ p\ne\one.
$$

\sevafig{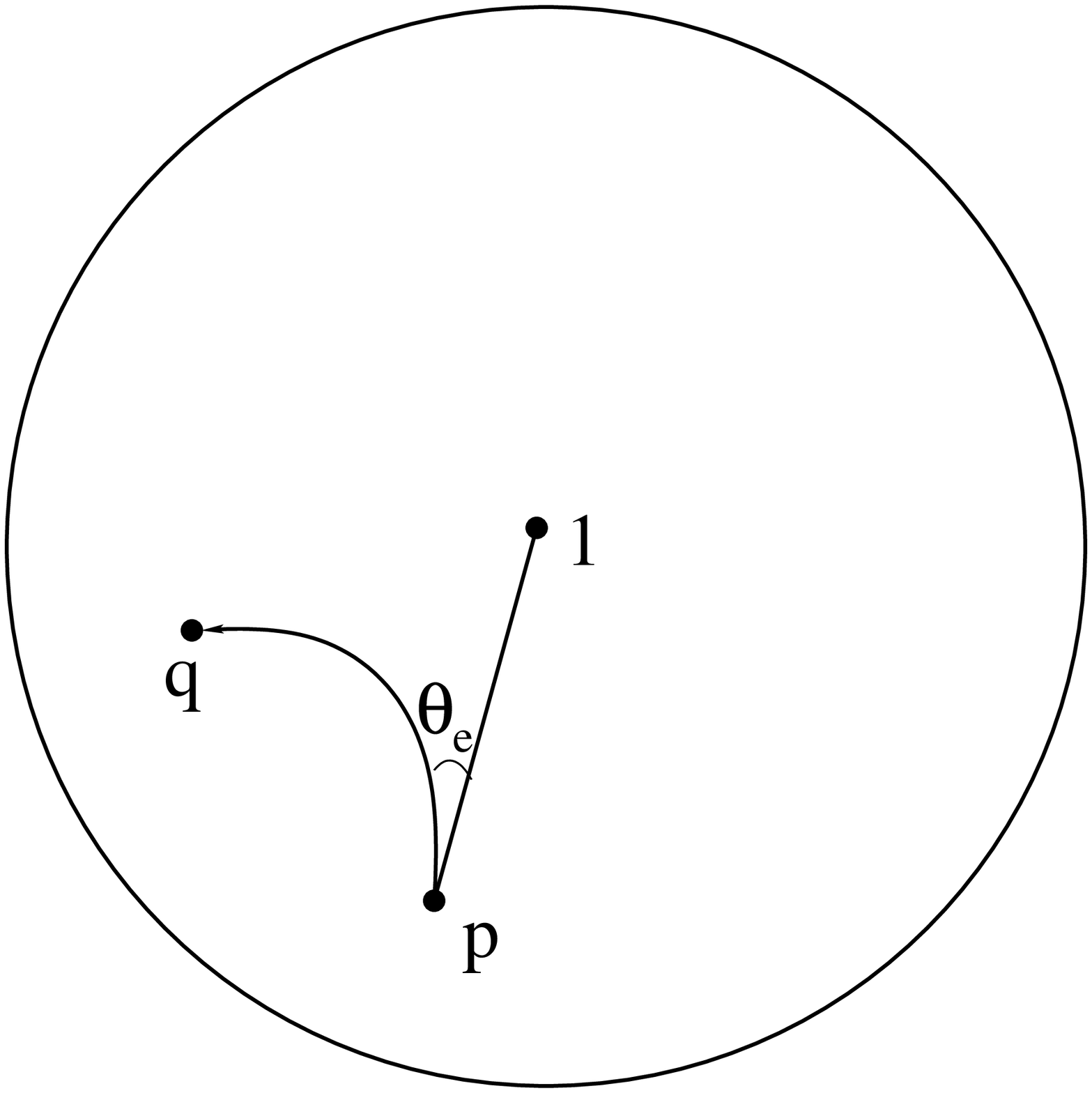}{70mm}{0}

In this case we consider the angle $\theta_e$ between the geodesic
$l(p,q)$ in the Poincar\'e metric on $D^2$, with the geodesic
$l(p,\one)$ (the last geodesic should be a diametral line), counted
from $l(p,\one)$ to $l(p,q)$ counterclockwise. The angle $\theta_e$ is
defined modulo~$2\pi$. By definition, $\varphi_e=d\theta_e$; it is a
well-defined $1$-form on the space $D_{\one,2,0}$, or on
$D_{\one,1,1}$, see Figure 1.

\item[C2)] $e=(p,q)$, $p=\one$.

\sevafig{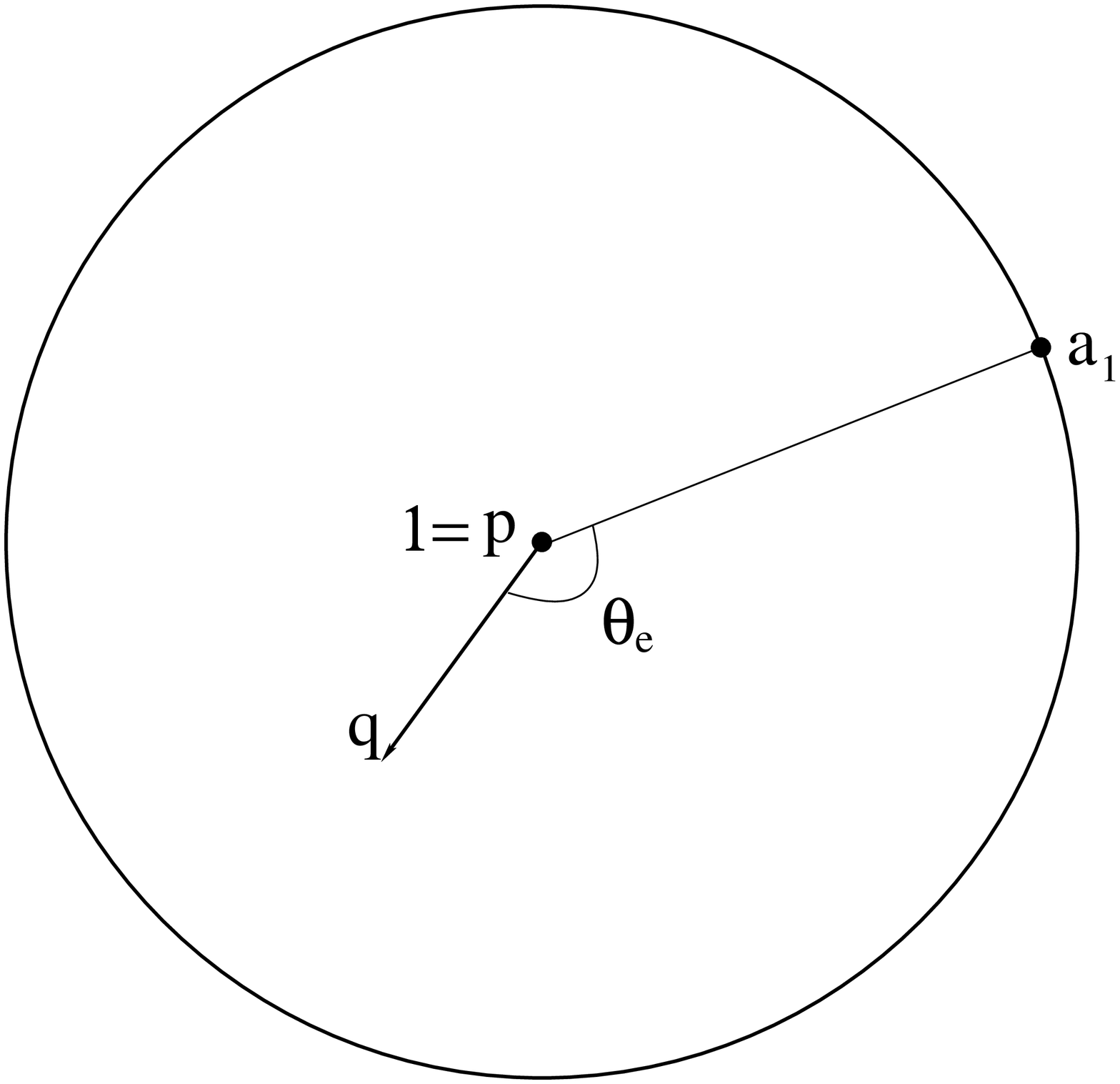}{60mm}{0}

In this case $\theta_e$ is the angle between the line $(\one,q)$ and
the line $(\one,a_1)$ where $a_1=\{\bar1\}$, the first vertex of the
second type on~$S^1$. 
 The $1$-form $\varphi_e=d\theta_e$ is well-defined $1$-form on
$D_{\one,1,1}$ or $D_{\one,0,2}$. We denote by $\varphi_e$ also its
pull-back on $D_{\one,n,m}$, see Figure 2.
\end{list}

\subsubsection{Polydifferential operators, corresponded to the graphs}

Now we want to define some operators
$$
\Lambda^nT^\ndot_\poly(\R^d)\otimes
C_\ndot(A,A)\to\Omega^\ndot(\R^d)[-n].
$$
Let $\gamma_1,\dots,\gamma_n$ be polyvector fields on $\R^d$, and let
$a_1,\dots,a_m$ be $m$ functions, $a_i\in A=C^\infty(\R^d)$. Now let
$\Gamma$ be an admissible graph, such that $\#\Star(i)=\deg\gamma_i+1$
for a vertex $i$ of the first type, $\#\Star(\one)=l$. Then we are
going to construct an outcoming  $l$-differential form
$\Omega^\Gamma_l$. To do that, we define
$\Omega^\Gamma_l(\partial_{\alpha_1}\wedge\dots\wedge\partial_{\alpha_l})$,
$i_1,\dots,i_l\in\{1,\dots,d\}$. As in Section~2.1.2,
$$
\Omega^\Gamma_l(\partial_{\alpha_1}\wedge\dots\wedge\partial_{\alpha_l})=\sum_{I\colon
E_\Gamma\setminus\Star(\one)\to\{1,\dots,d\}}\Omega^{\alpha_1,\dots,\alpha_l}_I
$$
We can extend the map
$I\colon
E_\Gamma\setminus\Star(\one)\to\{1,\dots,d\}$ to a map 
$\hat I\colon E_\Gamma\to\{1,\dots,d\}$, where if $e=(\one,{*})$
and $e$ has the label $e^s_\one$ in the graph~$\Gamma$, we set
$I(e)=\alpha_s$. Then the functions~$\Psi_v$, $v\ne\one$, are defined
as in Section~2.1.2. Next,
$$
\tilde\Psi_v=\left(\prod_{e\in E_\Gamma,\ e=({*},v)}\partial_{\hat
I(e)}\right)\Psi_v,\quad v\ne\one,
$$
and then
$$
\Omega^{\alpha_1,\dots,\alpha_l}_I=\prod_{v\in
V_\Gamma\setminus\one}\tilde\Psi_v.
$$
Finally,
$$
\Omega^\Gamma_l(\partial_{\alpha_1}\wedge\dots\wedge\partial_{\alpha_l})=\sum_I\Omega^{\alpha_1,\dots,\alpha_l}_I.
$$

\subsubsection{The grading}

Let us find a relation on $n,m,\#E_\Gamma$ and~$l$. We want to define
a map 
$$
\Lambda^nT^\ndot_\poly\otimes C_{m-1}(A,A)\to\Omega^\ndot[-n].
$$
Therefore, 
\begin{equation}
\deg\gamma_1+\dots+\deg\gamma_n+(-m+1)-n=-l.
\end{equation}
On the other hand
\begin{equation}
\#E_\Gamma=\deg\gamma_1+\dots+\deg\gamma_n+n+l.
\end{equation}
We have from (26) and (27):
\begin{equation}
\#E_\Gamma-n-l-m+1-n=-l
\end{equation}
or
\begin{equation}
\#E_\Gamma=2n+m-1.
\end{equation}

Therefore, for the grading condition a graph with $n$ vertices of the
first type and $m$ vertices of the second type should have $2n+m-1$
edges. According to our discussion in Section~2.1.4, the configuration
space for the construction of~$\hat\U$ should have dimension
$2n+m-1$. The space $D_{\one,n,m}$ has precisely this dimension.

\subsubsection{Weight $W_\Gamma$}

Now we consider a graph $\Gamma$ as placed on the space
$\Disk_{\one,n,m}$ such that the vertex $\one$ of $\Gamma$ is the
center $\one$ of the disk~$D^2$, the points $\{p_1,\dots,p_n\}$ of
$\Disk_{\one,n,m}$ are the vertices $\{1,\dots,n\}$ of $\Gamma$ of the
first type, and the points $\{q_1,\dots,q_m\}$ of $\Disk_{\one,n,m}$
are the vertices $\{\bar1,\dots,\bar m\}$ of $\Gamma$ of the second
type. Then each edge $e$ of~$\Gamma$ defines a map
\begin{align*}
&pr_e\colon\Disk_{\one,n,m}\to\Disk_{\one,2,0}\quad\text{or}\\
&pr_e\colon\Disk_{\one,n,m}\to\Disk_{\one,1,1}
\end{align*}
We denote by the same symbol $\varphi_e$ the pull-back
$\pr^*_e\varphi_e$ of the $1$-form, constructed in Section~2.2.3. The
$1$-form $\varphi_e$ is $G$-invariant ($G$ is the rotation group), and
we can consider $\varphi_e$ as a $1$-form on the space $\D_{\one,n,m}$.

\begin{defin}[the weight $W_\Gamma$]
\begin{equation}
W_\Gamma=\frac1{(\#\Star(\one))!}\prod_{k=1}^n\frac1{(\#\Star(k))!}\cdot\frac1{(2\pi)^{2n+m-1}}\cdot\int_{D^+_{\one,n,m}}\bigwedge_{e\in
E_\Gamma}\varphi_e
\end{equation}
where $D^+_{\one,n,m}$ is the connected component of $D_{\one,n,m}$
for which the points $(q_1,\dots,q_m)\in S^1$ define the right cyclic
order on~$S^1$.
\end{defin}

\subsubsection{}
\begin{theorem}
Let $G_{\one,n,m}$ be the set of the admissible graphs in the sence of
Section~\emph{2.2.2} with $2n+m-1$ edges. Define the map
$$
\hat\U_n\colon\Lambda^nT^\ndot_\poly(\R^d)\otimes
C_{m-1}(A,A)\to\Omega^\ndot(\R^d)[-n]
$$
by the formula
\begin{gather}
\hat\U_n=\sum_{\Gamma\in G_{\one,n,m}}W_\Gamma\times\Omega^\Gamma_l,\\
l=n+m-1-\sum^n_{i=1}\deg\gamma_i.
\end{gather}
Then the maps $\{\hat\U_n\}$ are the Taylor components of an
$L_\infty$-morphism of the $L_\infty$-modules over
$T^\ndot_\poly(\R^d)$\emph:
$$
\hat\U\colon C_\ndot(A,A)\to\Omega^\ndot(\R^d).
$$
The $L_\infty$-module over $T^\ndot_\poly(\R^d)$ structure on
$C_\ndot(A,A)$ is defined by \emph{(19)} where $\{\U_k\}$ are the Taylor
components of the Kontsevich formality morphism \emph(with the harmonic
angle function, as in Section~\emph{2.1)}.
\end{theorem}

\begin{proof}
We need to prove (21) for $\{\hat U_n\}$ defined as above. The l.h.s.\
of~(21) has a form
$$
\sum_{\Gamma\in
G'_{\one,n,m}}C_\Gamma\cdot\Omega^\Gamma_l(\gamma_1\wedge\dots\wedge\gamma_n\otimes(a_1\otimes\dots\otimes
a_m))
$$
where $G'_{\one,n,m}$ is the set of admissible graphs with $n$
vertices of the first type, $m$ vertices of the second type, and
$2n+m-2$ edges, and $l=n+m-2-\sum_{i=1}^n\deg\gamma_i$. The numbers
$C_\Gamma$ are linear-quadratic expressions in the
weights~$W_{\Gamma_1}$,defined in Section~2.2.6, and the Kontsevich
weights~$W_{\Gamma_2}$, defined in Section~2.1.6. We want to prove
that $C_\Gamma=0$ for any $\Gamma\in G'_{\one,n,m}$. The idea arise
to~\cite{K}: we identify $C_\Gamma$ with the integral over the
boundary $\partial\bar D_{\one,n,m}$ of the differential form of
degree $2n+m-2$, namely, with $\bigwedge_{e\in
E_\Gamma}d\varphi_e$. We have by the Stokes formula:
\begin{equation}
\int_{\partial\bar D_{\one,n,m}}\bigwedge_{e\in
E_\Gamma}d\varphi_e=\int_{        \bar D_{\one,n,m}}d\left(\bigwedge_{e\in
E_\Gamma}d\varphi_e\right)=0.
\end{equation}
On the other hand, only the boundary strata of codimension~$1$
contributes in the l.h.s.\ of~(33). We have:
\begin{equation}
\int_{\partial\bar D_{\one,n,m}}\bigwedge_{e\in
E_\Gamma}d\varphi_e=\int_{\partial_{\text{S1)}\bar D_{\one,n,m}}}\bigwedge_{e\in
E_\Gamma}d\varphi_e+\int_{\partial_{\text{S2)}\bar D_{\one,n,m}}}\bigwedge_{e\in
E_\Gamma}d\varphi_e+\int_{\partial_{\text{S3)}\bar D_{\one,n,m}}}\bigwedge_{e\in
E_\Gamma}d\varphi_e
\end{equation}
where S1), S2), and S3) are the three types of the boundary strata of
codimension~$1$, listed in Section~2.2.1. The idea of the proof is to
identify these summands with the summands contributed to $C_\Gamma$
from~(21). We consider the cases S1--S3) separately.

\subsubsubsection{The case S1).}
The boundary stratum is $C_k\times D_{\one,n-k+1,m}$. It is clear that
the integral over it factorizes in the product of two integrals: the
integral over $C_k$ and the integral over $D_{\one,n-k+1,m}$. It is
proved in~\cite{K}, Section~6, that the integral over~$C_k$ does not
vanish only for $k=2$. The situation is as follows: two points move
close to each other and far from~$\one$. This corresponds to the first
summand of~(21), containing the Schouten--Nijenhuis bracket of
polyvector fields.

\subsubsubsection{The case S2).}
Some points $p_{i_1},\dots,p_{i_k}$ move close to each other and close
to
$\one\in D^2$. 
The
boundary stratum is $D_k\times D_{\one,n-k,m}$. 
It is clear that the integral factorizes into the product of two
integrals, of an integral over~$D_k$ and of an integral over
$D_{\one,n-k,m}$ (see, in particular, the following two
subsections). Then, by Theorem 6.6.1 in~[K], we have only the following
possibilities:

\begin{list}{}{}
\item[S2.1)] $k=1$ and there is an edge from $\one$ to~$p_{i_1}$;
\item[S2.2)] $k=1$ and there are no edges connecting $p_{i_1}$
and~$\one$.
\end{list}

We analyze these strata separately below. It will be clear from our
discussion how to apply Theorem~6.6.1 in~[K] to prove that other
possible strata give zero contribution.

\subsubsubsubsection{The case S2.1).}

\sevafig{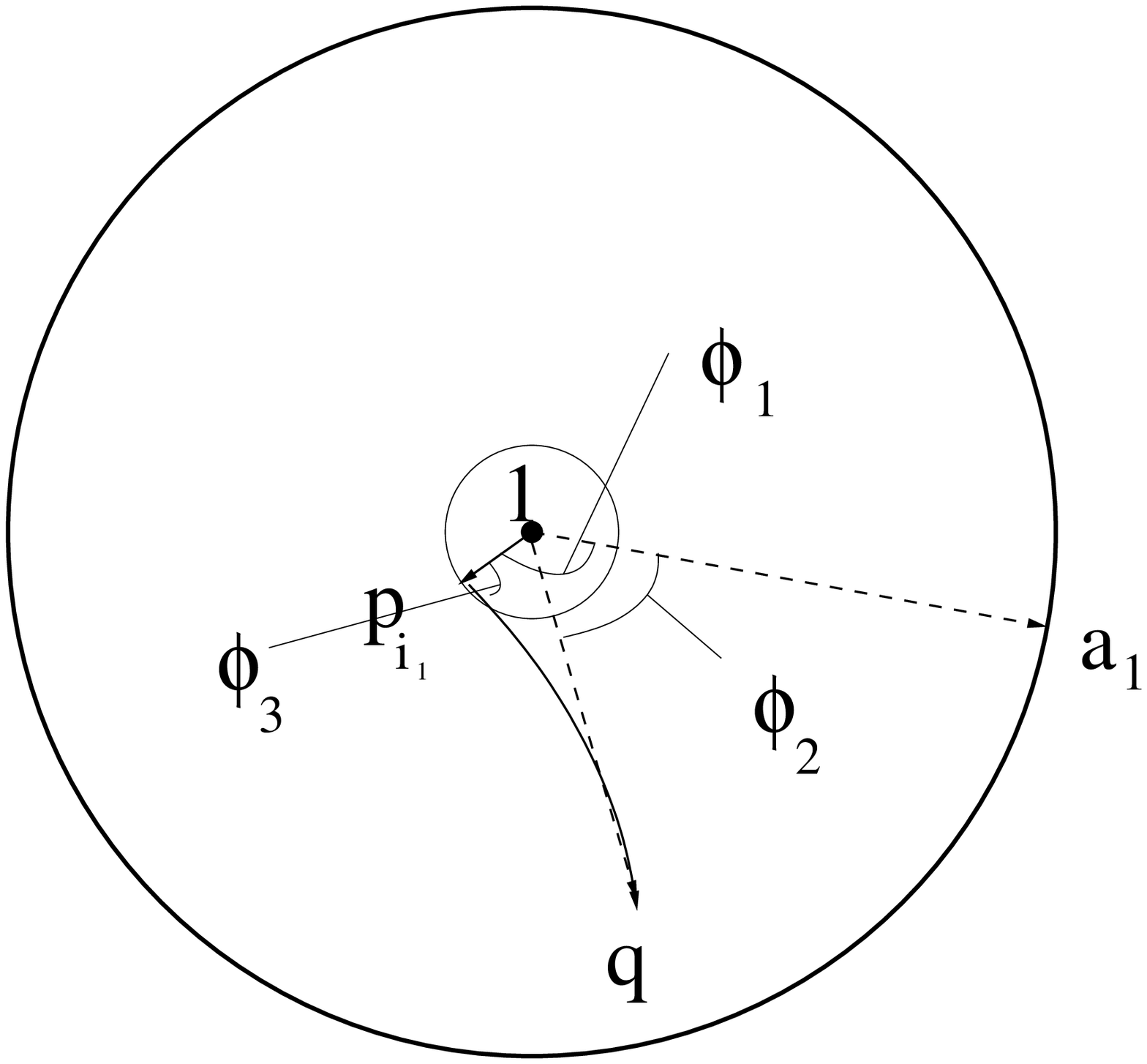}{90mm}{0}

Let $(p_{i_1},q)$ be an edge starting at~$p_{i_1}$, and let $\Phi_1$
be the angle between $(\one,p_{i_1})$ and $(\one,a_1)$ (this is the
angle $\theta_e$ for $e=(\one,p_{i_1})$), $\Phi_2$ be the angle
between $(\one,q)$ and $(\one,a_1)$, and $\Phi_3$ be the angle
between $(p_{i_1},q)$ and $(p_{i_1},\one)$, see Figure 3. The integral over $D_1$ should be the
integral over~$D_1$ of the canonical $1$-form $d\varphi$ on~$D_1$. By
our definition, $\theta_e$~for $e=(p_{i_1},q)$ is the
angle~$\Phi_3$. The $2$-form $d\varphi_{(p_{i_1},q)}\wedge
d\varphi_{(\one,p_{i_1})}$ is equal to $d\Phi_3\wedge d\Phi_1$. In the
limit $p_{i_1}\to\one$ the angle $\Phi_3\to\pi-(\Phi_1-\Phi_2)$ (see
Fig.~3). Therefore, we have two-form $d(\pi-\Phi_1+\Phi_2)\wedge
d\Phi_1=d\Phi_2\wedge d\Phi_1$. The angle $\Phi_2$ does not depend on
the position of~$p_{i_1}$. The ``new'' edge will be $(\one,q)$, and
its angle is precisely~$\Phi_2$. The integral over $D_1$ is $\int
d\Phi_1$.

In the general case, there are some edges starting at $p$, but no edge ends at $p$: in this case the angle of this edge would be equal to 0.

These terms correspond to a summand of $d\circ
i_{\gamma_i}(\hat\U_{k+1}(\gamma_{j_1}\wedge\dots\wedge\gamma_{j_{k+1}}\otimes\omega))$
(see~(21)) and any they contribute to the third summand in~(21). See
Section~2.2.7.2.3 for a further discussion.

\subsubsubsubsection{The case S2.2).}
We denote $p_{i_1}$ by $p$. There no edges connecting $p$
and~$\one$. See Figure~4.

\sevafig{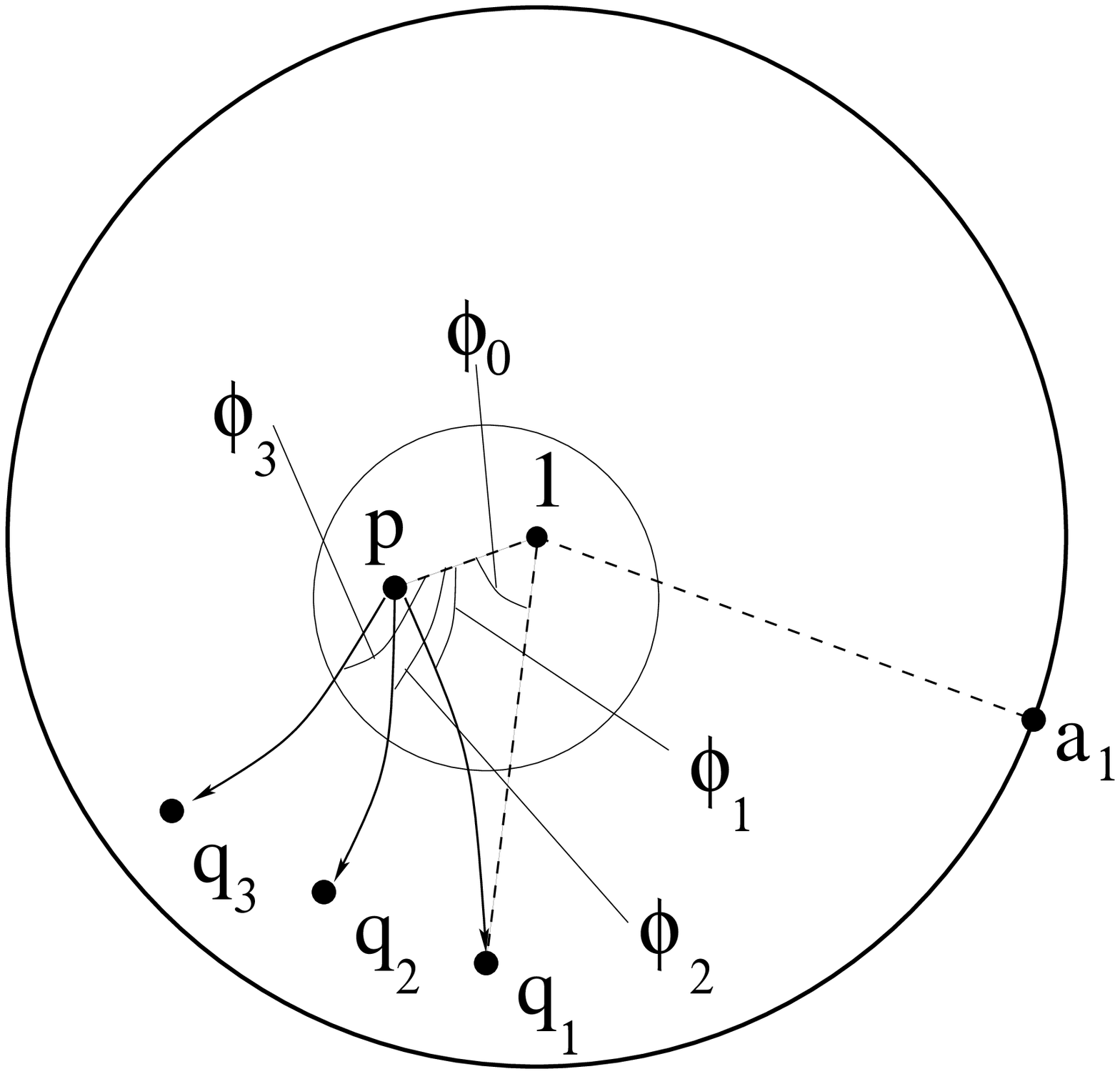}{90mm}{0}

In the limit when $p$ tends to~$\one$ and the points
$q_1,q_2,q_3\dots$ (end-points  of the edges starting from~$p$) are
far from~$\one$, the wedge product $d\phi_1\wedge d\phi_2\wedge
d\phi_3\wedge\dots$ tends to $\pm d\phi_0\wedge d(\phi_2-\phi_1)\wedge
d(\phi_3-\phi_1)$. The integral $\int_{S^1}d\phi_0$ is the integral
over~$D_1$. The remaining part of the wedge product,
$d(\phi_2-\phi_1)\wedge d(\phi_3-\phi_1)\wedge\dots$ contributes to
the integral over $D_{\one, n-1,m}$.

In the general case, if there exists any edge ending at~$p$, an edge
$l=(s,p)$, and $s$ is far from~$\one$, the angle $\varphi_e$ tends
to~$0$ when $p$ tends to~$\one$. Therefore, no edge ends at~$p$.

These terms correspond to a summand of $i_{\gamma_i}\circ
d(\hat\U_{k+1}(\gamma_{j_1}\wedge\dots\wedge\gamma_{k+1}\otimes\omega))$
and contribute to the third summand in~(21). See Section 2.2.7.2.3 for
a further discussion.

\subsubsubsubsection{}
We claim that the total contribution of the two cases S2.1) and S2.2)
is equal exactly to $L_{\gamma_i}\tilde\omega=d\circ
i_{\gamma_i}\tilde\omega\pm i_{\gamma_i}\circ d\tilde\omega$, where
$\tilde\omega=\hat\U_{k+1}(\gamma_{j_1}\wedge\dots\wedge\gamma_{j_{k+1}}\otimes\omega)$.
Here we consider two examples.

\begin{nexample}
$\gamma_i=v$ is a vector field, $\tilde\omega$ is a function. Then
$d\circ i_v\tilde\omega=0$, and we have only the summand $i_v\circ
d\tilde\omega$, which corresponds exactly exactly to the case~S2.2).
\end{nexample}

\begin{nexample}
Again, $\gamma_i=v$ is a vector  field, $v=\sum_iv_i\partial_i$, but
$\tilde\omega$ is a $1$-form $\tilde\omega=\sum_if_idx_i$. Then
\begin{equation}
d\circ
i_v\tilde\omega=d\left(\sum_iv_if_i\right)=\sum_{i,j}f_i\frac{\partial
v_i}{\partial x_j}dx_i+\sum_{i,j}v_i\frac{\partial f_i}{\partial
x_j}dx_j
\tag{$*$}
\end{equation}
On the other hand,
\begin{equation}
i_vd\omega=i_v\left(\sum_{i,j}\frac{\partial f_i}{\partial
x_j}dx_j\wedge dx_i\right)=\sum_{i,j}\frac{\partial f_i}{\partial
x_j}\cdot v_jdx_i-\sum_{i,j}\frac{\partial f_i}{\partial x_j}
v_idx_j
\tag{$**$}
\end{equation}
We see that in $L_v\tilde\omega=d\circ
i_v\tilde\omega+i_v\circ d\tilde\omega$ the second summand in the
r.h.s.\ of~($*$) cancells the second summand in the r.h.s of~($**$),
whereas the first summand in the r.h.s.\ of~($*$) corresponds to the
case~S2.1), and the first summand in the r.h.s.\ of~($**$) corresponds
to the case~S2.2).
\end{nexample}

In the general case the situation is analogous to this example.

\subsubsubsection{The case S3).}
Some points $p_{i_1},\dots,p_{i_k}$ of the first type and some points
$q_{j_1},\dots,q_{j_l}$ of the second type move close to each other
and close to $S^1=\partial\bar{D^2}$, $2k+l\ge2$. The boundary stratum
is $C_{k,l}\times D_{\one,n-k, m-l+1}$. The integral over $C_{k,l}$ is
exactly the Kontsevich integral: the angle with~$\one$ becomes the
angle with $\infty$ for $C_{k,l}$. We obtain the second summand in~(21).

Theorem 2.2.7 is proven.
\end{proof}

\subsubsection{}
It is clear that
$$
\hat\U_0(a_1\otimes\dots\otimes a_m)=\mu(a_1\otimes\dots\otimes
a_m)=\frac1{(m-1)!}a_1\,da_2\wedge\dots\wedge da_m.
$$

\section{Applications to traces, cup-products, and to the Duflo formula}

\subsection{}
First of all, we recall what the (extended) Duflo formula is. Let $\g$ be a finite-dimensional Lie algebra, we associate with $\g$ two $\g$-modules: the symmetric algebra $S^\ndot(\g)$, and the universal enveloping algebra $U(\g)$. These two $\g$-modules are isomorphic due to the Poincar\'e--Birkhoff--Witt theorem, the isomorphism $\varphi_\PBW\colon S^\ndot(\g)\to U(\g)$ is defined as
$$
\varphi_\PBW(g_1\cdot\dots\cdot g_k)=\frac1{k!}\sum_{\sigma\in\Sigma_k}g_{\sigma(1)}\otimes\dots\otimes g_{\sigma(k)}.
$$
As algebras, $S^\ndot(\g)$ and $U(\g)$ are not isomorphic. The Duflo
theorem states that the algebras of $\g$-invariants $[S^\ndot(\g)]^\g$
and $[U(\g)]^\g$ are (canonically) isomorphic. Let us recall the
construction of this isomorphism.

For any $k\ge1$, let $\Tr_k$ be a canonical invariant element
$\Tr_k\in[S^k(g^*)]^\g$, defined as the symmetrization of the map
$g\mapsto\Tr|_\g\ad^kg$ ($g\in\g$). We consider elements of
$S^k(\g^*)$ as differential operators of $k$-th order on
$S^\ndot(\g)$ with constant coefficients.  We define the
map $\varphi_\strange\colon S^\ndot(\g)\to S^\ndot(\g)$ as
$$
\varphi_\strange=\exp\left(\sum_{k\ge1}\alpha_{2k}\Tr_{2k}\right)
$$
where the rational numbers $\alpha_{2k}$ are defined as
$$
\sum_{k\ge1}\alpha_{2k}x^{2k}=\frac12\log\frac{e^{\frac
x2}-e^{-\frac x2}}x.
$$
It is clear that $\varphi_\strange\colon S^\ndot(\g)\to S^\ndot(\g)$
is a map of $\g$-modules, because $\Tr_{2k}\in[S^{2k}(\g^*)]^\g$.

\subsubsection{}

\begin{theorem}[Duflo]
For any finite-dimensional Lie algebra $\g$ the map
$$
\varphi^*_D=[\varphi_\PBW\circ\varphi_\strange]\colon[S^\ndot(\g)]^\g\to[U(\g)]^\g
$$
is an isomorphism of algebras.
\end{theorem}

\subsubsection{}

\begin{theorem}[Kontsevich~\cite{K}]
For any finite-dimensional Lie algebra $\g$ the map $\varphi_D$ of
$\g$-modules $\varphi_D\colon S^\ndot(\g)\to U(\g)$ defines the map of
algebras
$$
\varphi_D^*\colon H^\ndot(\g;S^\ndot(\g))\to H^\ndot(\g;U(\g)).
$$
\end{theorem}

\subsubsection{}
Now consider \emph{homology} $H_\ndot(\g;S(\g))$ and
$H_\ndot(\g;U(\g))$. They have a structure of modules over the
algebras $H^\ndot(\g,S(\g))$ and $H^\ndot(\g;U(\g))$,
correspondingly. In particular, for zero (co)homology
$A^\coinv=A/\g\cdot A$ has a structure of a module over the algebra
$A^\inv$ for any $\g$-module $A$ which is an associative algebra such
that the multiplication $A\otimes A\to A$ is a map of $\g$-modules.

We want to construct a map of modules over $H^\ndot(\g,\dots)$
$$
\varphi_{D{*}}\colon H_\ndot(\g;S(\g))\to H_\ndot(\g;U(\g))
$$
such that it is compatible with the map of algebras
$$
\varphi_D^*\colon H^\ndot(\g;S(\g))\to H^\ndot(\g,U(\g)).
$$
It means that for any $\omega\in H^\ndot(\g;S^\ndot(\g))$ and $\eta\in
H_\ndot(\g,S^\ndot(\g))$ one has
\begin{equation}
\varphi_{D{*}}(\omega\cdot\eta)=\varphi_D^*(\omega)\cdot\varphi_{D{*}}(\eta).
\end{equation}

\subsubsection{}
\begin{conjecture}
For any finite-dimensional Lie algebra $\g$, the map $\varphi_{D{*}}$,
induced by the map $\varphi_D=\varphi_\PBW\circ\varphi_\strange$ of
$\g$-modules, satisfies~\emph{(35)}.
\end{conjecture}

This conjecture is a typical application of a conjecture on
cup-products on tangent cohomology, Conjecture 3.5.3.1 below. For a
semisimple (or, more generally, unimodular) Lie algebra $\g$
Conjecture 3.1.4 follows from Theorem 3.1.2: for example, the
top-cohomology $H^{top}(\g;A)\simeq H_0(\g;A)$, and this is is an isomorphism of
$H^0(\g;A)$-modules for any $A$ with compatible structures of an
associative algebra and a $\g$-module. Analogously for higher
(co)homology.

For a general Lie algebra $\g$, $H^{top}(\g;M)\simeq H_0(\g;\Tr\otimes
M)$ where $\Tr\colon\g\to\C$ is a one-dimensional $\g$-module, the
trace of the adjoint action. Therefore, for a general Lie algebra $\g$, Theorem~3.1.2 does not imply
Conjecture~3.1.4.

\subsection{$L_\infty$-quasi-isomorphisms of $L_\infty$-modules, the
tangent complexes, the tangent map}

We say that an $L_\infty$-map of two $L_\infty$-modules
$M^\ndot,N^\ndot$ (see Section~1.5.4) is an
\emph{$L_\infty$-quasi-isomorphism of modules} if $\varphi_0\colon
M^\ndot\to N^\ndot$ (which a priori is a map of complexes) is a
quasi-isomorphism of complexes.

Let $M^\ndot$ be an $L_\infty$-module over a dg Lie algebra~$\g$, and
let $\pi\in\g^1$ satisfies the Maurer--Cartan equation
$d\pi+\frac12[\pi,\pi]=0$. We define the \emph{tangent complex} $T_\pi
M^\ndot$ as the graded space $M^\ndot$ with the differential $d_\pi$
defined as follows:
\begin{equation}
d_\pi(m)=\phi_0(m)+\phi_1(\pi,m)+\frac12\phi_2(\pi,\pi,m)+\dots+\frac1{n!}\phi_n(\pi,\pi,\dots,\pi,m)+\dots
\end{equation}
where $\phi_k\colon\Lambda^k\g\otimes M^\ndot\to M^\ndot[1-k]$ are the
Taylor components of the $L_\infty$-module structure on~$M^\ndot$ (see
Section~1.5.3).

\begin{lemma}
\leavevmode

\begin{enumerate}
\item $d^2_\pi=0$ if $d\pi+\frac12[\pi,\pi]=0$\emph;
\item the formula
\begin{equation}
(T_\pi\varphi)(a)=\sum_{k=0}^\infty\varphi_{k+1}(\pi,\pi,\dots,\pi,a)/k!
\end{equation}
defines a map of the tangent \emph{complexes}
$$
T_\pi\varphi\colon T_\pi M^\ndot\to T_\pi N^\ndot;
$$
\item let $\varphi$ be an $L_\infty$-quasi-isomorphism of
modules. Then $T_\pi\varphi$ is a quasi-isomorphism of the complexes
provided the condition that $\pi$ is sufficiently small. In
particular, if we replace $\pi$ on $\hbar\pi$, where $\hbar$ is a
formal parameter, $T_{\hbar\pi}\varphi$ is a quasi-isomorphism.
\end{enumerate}
\end{lemma}

\begin{proof}
It is straightforward.
\end{proof}

\subsection{}
\relax   From now on, we will work with formal power series in a formal
parameter~$\hbar$.

\begin{lemma}
\leavevmode

\begin{enumerate}
\item Let $A=C^\infty(\R^d)[[\hbar]]$,
$\pi\in\Lambda^2T_\poly(\R^d)$, $[\pi,\pi]=0$. Then the map
$$
T_{\hbar\pi}\hat\U\colon T_{\hbar\pi}C_\ndot(A,A)\to
T_{\hbar\pi}\Omega^\ndot(\R^d)[[\hbar]]
$$
is a quasi-isomorphism of the complexes\emph;
\item 
$$
T_{\hbar\pi}C_\ndot(A,A)=C_\ndot(A_*,A_*)
$$
where $A_*$ is the deformed algebra~$A$ with the Kontsevich
star-product\emph;
\item
$$
T_{\hbar\pi}\Omega^\ndot(\R^d)[[\hbar]]=(\Omega^\ndot(\R^d)[[\hbar]],
d=l_{\hbar\pi}).
$$
\end{enumerate}
\end{lemma}

\begin{proof} (i) follows from Lemma 3.1, because $\hat\U_0=\mu$ is a
quasi-isomorphism by Theorem~1.3.1.

(ii) and (iii) follow from the definitions.
\end{proof}

In this way, we obtained a quasi-isomorphism of complexes:
$$
T_{\hbar\pi}\hat\U\colon
C_\ndot(A_*,A_*)\to(\Omega^\ndot(\R^d)[[\hbar]],
d=L_{\hbar\pi}).
$$

\subsection{}
\begin{lemma}
\leavevmode

\begin{enumerate}
\item
$$
H_0(C_\ndot(A_*,A_*))=A_*/[A_*,A_*];
$$
\item
$$
H_0(\Omega^\ndot(\R^d)[[\hbar]],d=L_{\hbar\pi})=A/\{A,A\},
$$
where $\{A,A\}$ is the commutant of the Poisson algebra.
\end{enumerate}
\end{lemma}

\begin{proof}
(i) In $C_\ndot(A_*,A_*)$ we have:
$$
d(a_0\otimes a_1)=a_0*a_1-a_1*a_0.
$$

(ii) By definition,
$\Im(d\colon\Omega^1\to\Omega^0)=L_{\hbar\pi}(\Omega^1)$. We have:
$$
L_{\hbar\pi}(fdg)=(d_\DR\circ i_{\hbar\pi}+i_{\hbar\pi}\circ
d_\DR)(fdg)=i_{\hbar\pi}\circ d_\DR(fdg)=\hbar\pi(df\wedge
dg)=\hbar\{f,g\}.
$$
\end{proof}

\subsubsection{}
\begin{theorem}
The map $T_{\hbar\pi}\hat\U$ gives an isomorphism $A_*/[A_*,A_*]\simto
A/\{A,A\}$ for any Poisson bivector~$\pi$.\qed
\end{theorem}

We return to this map for a linear Poisson structure $\pi$ in Section 3.6 after a
conjecture on the compatibility with cup-products on the level of
cohomology.

\subsection{Cup-products on the tangent cohomology}

\subsubsection{}

\begin{defin}
For a dg Lie algebra~$\g$, and a solution $\pi$ of the Maurer--Cartan
equation $d\pi+\frac12[\pi,\pi]=0$, the tangent complex $T_\pi\g$ is
$\g[1](d_{\g}+\ad\pi)$.
\end{defin}

\begin{lemma}[\cite{K},\dots]
\leavevmode\par
\begin{enumerate}
\item Let $\U\colon\g^\ndot_1\to\g^\ndot_2$ be an $L_\infty$-map
between two dg Lie algebras, $\pi\in\g^1_1$ be a solution of the
Maurer--Cartan equation. Then
\begin{equation}
\tilde\pi:=\U_1(\pi)+\frac12\U_2(\pi,\pi)+\dots+\frac1{n!}\U_n(\pi,\dots,\pi)+\dots
\end{equation}
is a solution of Maurer--Cartan equation in~$\g^\ndot_2$ \emph(it is
clear that $\tilde\pi\in\g^1_2$\emph{);}
\item the map
$$
T_\pi\U\colon T_\pi\g_1\to T_{\tilde\pi}\g_2,
$$
defined as
\begin{multline}
T_\pi\U(x)=\U_1(x)+\U_2(x,\pi)+\\
+\frac12\U_3(x,\pi,\pi)+\dots+\frac1{(n-1)!}\U_n(x,\pi,\pi,\dots,\pi)+\dots
\end{multline}
is a map of complexes.
\end{enumerate}
\end{lemma}

\subsubsection{}
In the case of the formality morphism $\U\colon
T^\ndot_\poly(\R^d)\to\D^\ndot_\poly(\R^d)$ both dg algebras have an
extra property: the dg spaces $T_\pi T^\ndot_\poly(\R^d)$ and
$T_{\tilde\pi}\D^\ndot_\poly(\R^d)$ have structures of associative
algebras. In the case $T_\pi T^\ndot_\poly(\R^d)$ it is just the wedge
product of polyvector fields; in the case
$T_{\tilde\pi}\D^\ndot_\poly(\R^d)$ it is the usual cup-product in the
deformed algebra:
\begin{equation}
(\Psi_1\cup\Psi_2)(a_1\otimes\dots\otimes
a_{k+l})=\Psi_1(a_1\otimes\dots\otimes
a_k)*\Psi_2(a_{k+1}\otimes\dots\otimes a_{k+l}).
\end{equation}
The following remarkable result is proved in~\cite{K}, Section~8:

\begin{theorem}[Kontsevich]
The tangent map $T_\pi\U$ induces a map of the associative algebras on
the level of cohomology.
\end{theorem}

In other words, the map
$$
[T_\pi\U]\colon H^\ndot(T_\pi T^\ndot_\poly)\to
H^\ndot(T_{\tilde\pi}\D^\ndot_\poly)
$$
is a map of algebras.

Actually, Kontsevich considers a bit more general situation: a solution
$\pi$ of the Maurer--Cartan equation in $(\g^\ndot\otimes\m^\ndot)^1$
where $\m^\ndot$ is a commutative dg algebra. This generality is very
useful for applications, but we restrict ourselves by the case
$\m^\ndot=\C[0]$.

\subsubsection{} 
Now the tangent complex
$T_{\hbar\pi}(C_\ndot(A,A))=C_\ndot(A_*,A_*)$ has a module structure
over $C^\ndot(A_*,A_*)=T_{\hbar\pi}C^\ndot(A,A)$ (for any
algebra~$A_*$):

For $\Psi\colon A^{\otimes k}_*\to A_*$ and for $\omega=a_0\otimes
a_1\otimes\dots\otimes a_n$ ($k\le n$) we have the ``cup-product'':
\begin{equation}
\Psi(\omega)=(a_0*\Psi(a_1\otimes\dots\otimes a_k))\otimes
a_{k+1}\otimes\dots\otimes a_n.
\end{equation}

\begin{lemma}
\leavevmode\par
\begin{enumerate}
\item The ``cup-product'' described above
$$
\cup\colon C^\ndot(A_*,A_*)\otimes C_\ndot(A_*,A_*)\to C_\ndot(A_*,A_*)
$$
is a map of the complexes\emph;
\item
It endowes $C_\ndot(A_*,A_*)$ with an algebra structure over
$(C^\ndot(A_*,A_*))^\opp$ \emph(the algebra $A^\opp$ is the algebra
with the opposite multiplication\emph{: $a*_{A^\opp}b=b*_Aa$)}.
\end{enumerate}
\end{lemma}

\begin{proof}
It is straightforward.
\end{proof}

The map
$$
\cup\colon
T^\ndot_\poly(\R^d)\otimes\Omega^\ndot(\R^d)\to\Omega^\ndot(\R^d)
$$
(cup-product for $\Omega^\ndot$) is defined as the operator $i_{\gamma}$ of
the insertion of the polyvector field in the differential form. It is a
map of complexes:
$$
\cup\colon(T^\ndot_\poly(\R^d),d=\ad\pi)\otimes(\Omega^\ndot(\R^d),d=L_\pi)\to(\Omega^\ndot(\R^d),d=L_\pi).
$$
It is exactly the structure induced by (41) on the level of
cohomology.

\subsubsubsection{}

\begin{conjecture}
On the level of cohomology the tangent map
$$
T_{\hbar\pi}\hat\U\colon T_{\hbar\pi}C_\ndot(A,A)\to
T_{\hbar\pi}\Omega^\ndot(\R^d)[[h]]
$$
is a map of the modules, i.e.\ for $\omega\in
T_{\hbar\pi}C_\ndot(A,A)$ and for $\eta\in T_{\hbar\pi}T_\poly(\R^d)$
one has\emph:
\begin{equation}
[T_{\hbar\pi}\hat\U]([T_{\hbar\pi}\U]([\eta])\cup[\omega])=[\eta]\cup[T_{\hbar\pi}\hat\U]([\omega])
\end{equation}
\emph($[\dots]$ stands for the cohomological class of an element or
for the map induced on cohomology\emph).
\end{conjecture}

At the moment I don't know any proof of this Conjecture.


\subsection{}
Here we consider the case of a linear Poisson structure $\pi$
on~$\R^d$, $\pi=\sum c_{ij}^kx_k\partial_i\wedge\partial_j$. Recall,
that it follows from the equation $[\pi,\pi]=0$ that $\R^d\simeq\g^*$
for a Lie algebra~$\g$, and the numbers $\{c_{ij}^k\}$ are the
structure constants of this Lie algebra.

\subsubsection{}
First, let us describe the map $T_\pi\hat\U|_{A=C_0(A,A)}$ in this
case. We have $1$ point on the circle and some points inside the
disk~$D^2$. One can prove that all possible graphs are ``several
wheels,'' as it is shown in Fig.~5.

\sevafig{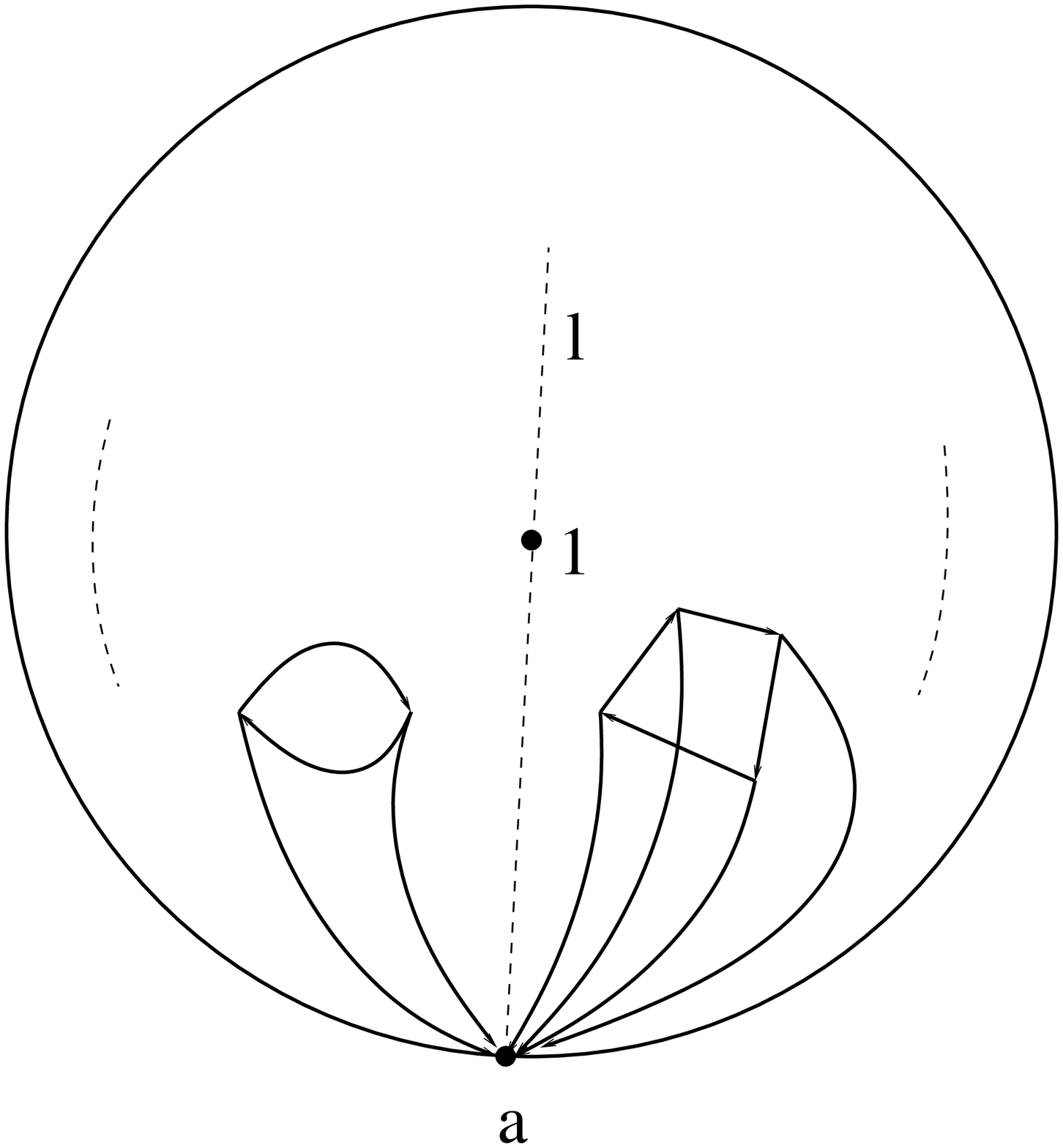}{80mm}{0}

\noindent
(There are no edges starting at $\one$ because we should obtain a
$0$-form, i.e.\ a function.) All the wheels have an even number of
vertices, because of the symmetry with respect to line~$l$ (see
Fig.~5). Let $w_{2k}$ be the weight of the single wheel with $2k$
vertices (in our sense, see Section~2.2.6).

We have:
\begin{equation}
(T_{t\pi}\hat\U)(a)=\exp\left(\sum_{k\ge1}\hbar^{2k}w_{2k}\Tr_{2k}\right)a
\end{equation}
(here $a\in S^\ndot(\g)=C_0(S^\ndot(\g),S^\ndot(\g))$). The operator
(43) satisfies the following properties:

\begin{enumerate}
\item it maps $[S(\g)_*,S(\g)_*]$ to $\{S(\g),S(\g)\}$;
\item if Conjecture 3.5.3.1 is true,
$$
[T_{\hbar\pi}\hat\U](\omega\cdot\eta)=\omega\cdot[T_{\hbar\pi}\hat\U](\eta)
$$
for $\omega\in[S(\g)]^\g$ and $\omega\in S(\g)_*/[S(\g)_*,S(\g)_*]$.
\end{enumerate}

The property (ii) follows from the ``vanishing of the
wheels''~\cite{Sh}, according to which $T_{\hbar\pi}\U=\Id$.

\begin{theorem}
If Conjecture \emph{3.5.3.1} is true, all the numbers $w_{2k}$,
$k\ge1$ are equal to zero.
\end{theorem}

\begin{proof}
The numbers $w_{2k}$, $k\ge1$ in (43) do not depend on the Lie
algebra~$\g$. Therefore, one can suppose that $\g$ is semisimple. We proved in
Section~3.1 that if $\g$ is semisimple, that for
$\omega\in[S^\ndot(\g)]^\g$ and $\eta\in S^\ndot(\g)/[\g,S^\ndot(\g)]$
one has:
$$
\varphi_{D{*}}(\omega\cdot\eta)=\varphi_D^*(\omega)\cdot\varphi_{D{*}}(\eta)
$$
where $\varphi_{D{*}}\colon[S^\ndot(\g)]^\coinv\to[U(\g)]^\coinv$ and
$\varphi_D^*\colon[S^\ndot(\g)]^\inv\to[U(\g)]^\inv$ are both induced
by the map $\varphi_D=\varphi_\PBW\circ\varphi_\strange$. Next,
consider the isomorphism $\Theta\colon S(\g)_*\to U(\g)$, defined as 
$$
\Theta(g_1*\dots*g_k)=g_1\otimes\dots*g_k.
$$
It is proven in~\cite{Sh} that
$\Theta=\varphi_\PBW\circ\varphi_\strange$.

We obtain that
\begin{equation}
\Id_*(\omega\cdot\eta)=\Id^*(\omega)\cdot \Id_*(\eta)
\end{equation}
for $\omega\in[S^\ndot(\g)_*]^\inv$, $\eta\in[S^\ndot(\g)]^\coinv$,
and $\Id_*\colon[S^\ndot(\g)_*]^\coinv\to[S^\ndot(\g)]^\coinv$, 
$\Id^*\colon[S^\ndot(\g)_*]^\inv\to[S^\ndot(\g)]^\inv$ are induced
by the identity map $\Id\colon S^\ndot(\g)_*\to S(\g)$. Suppose now that for
a map $T_*$ we have
\begin{equation}
T_*(\omega\cdot\eta)=\Id^*(\omega)\cdot T_*(\eta)
\end{equation}
($T_*=T_{\hbar\pi}(\hat\U)$ in our case). Then we have
$T_*=\Id_*$. Indeed, set $\eta=[1]$ and use the decompositions
\begin{align*}
U(\g)&=Z(U(\g))\oplus[U(\g),U(\g)],\\
S^\ndot(\g)&=[S^\ndot(\g)]^\inv\oplus\{S(\g),S(\g)\}
\end{align*}
which hold for a semisimple Lie algebra~$\g$.
\end{proof}

\subsubsection{}
\begin{corollary}
If Conjecture \emph{3.5.3.1} is true, then for any Lie algebra~$\g$
Conjecture~\emph{3.1.4} holds for $0$-\emph(co\emph)homology.
\end{corollary}

\subsubsection{}
Independently of Conjecture 3.5.3.1, the identity map $\Id\colon
S^\ndot(\g)_*\to S^\ndot(\g)$ is a map of $\g$-modules for any Lie
algebra~$\g$. Indeed, the map
$\varphi_D=\varphi_\PBW\circ\varphi_\strange\colon S^\ndot(\g)\to
U(\g)$ is a map of $\g$-modules, and the map $\varphi^{-1}_D\colon
U(\g)\to S^\ndot(\g)_*$ is an isomorphism of algebras (according
to~\cite{Sh}). Therefore, the composition
$\varphi^{-1}_D\circ\varphi_D=\Id$ is a map of $\g$-modules. In
particular, it defines a map of coinvariants
$\Id_*\colon[S^\ndot(\g)_*]^\coinv\to[S^\ndot(\g)]^\coinv$, i.e.\
$\Id$ maps $[S^\ndot(\g)_*,S^\ndot(\g)_*]$ to
$\{S^\ndot(\g),S^\ndot(\g)\}$. In other words, if we set
\begin{equation}
[f,g]_*=\sum_{k\ge0}\hbar^{2k+1}C_{2k+1}(f,g)
\end{equation}
then $C_{2k+1}=\{F_{2k+1},G_{2k+1}\}+\{F'_{2k+1},G'_{2k+1}\}+...$ for some $F_{2k+1},G_{2k+1},F'_{2k+1},G'_{2k+1},...\in
S^\ndot(\g)$. For example, $C_1(f,g)=\{f,g\}$.

\medskip

It would be very interesting to calculate the map
$T_{\hbar\pi}(\hat\U)$ for higher cohomology and to deduce
Conjecture~3.1.4 from Conjecture~3.5.3.1. The proof of
Conjecture~3.5.3.1 should be somewhat very close to the proof of the
Kontsevich theorem on cup-products on tangent cohomology~\cite{K},
Section~8.

\medskip

Finally, we have seen in this Section that the weels $w_{2k}$, $k\ge1$,
are conjecturally equal to zero, as well as the Kontsevich wheels~\cite{Sh}. It is very
interesting to compare the Kontsevich integrals of ``higher wheels``
and our integrals of them (a ``higher wheel'' $=$ a graph which appear
in $T_{\hbar\alpha}(\hat\U)(f)$ or $T_{\hbar\alpha}(\U)(f)$ for a non-linear~$\alpha$).
 
\section*{Acknowledgments}

The first proof of the Tsygan formality conjecture for chains, based
on operadic methods, appeared in the talk of Dima Tamarkin~\cite{T} at
the Mosh\'e Flato 2000 Conference. The author was inspired by this
talk and by further discussions with Boris Tsygan. Many discussions with Giovanni Felder on the first version of this paper were very useful, they helped me to find the omitted in the first version boundary stratum S2.2. I thank Seva Kordonsky for a quick and quality typing of this text. I am grateful to
IPDE grant 1999--2001 for the partial financial support and to the
ETH-Zentrum (Z\"urich) where the work was done for  hospitality and for
the remarkable atmosphere.

 \end{document}